\documentclass[10pt]{article}
\usepackage{amsmath,amsthm,amssymb,amsfonts,amscd}
\theoremstyle{plain}
\newtheorem{theorem}{Theorem}   
\newtheorem{proposition}[theorem]{Proposition}
\newtheorem{lemma}[theorem]{Lemma}
\newtheorem{corollary}[theorem]{Corollary}

\newtheorem{definition}[theorem]{Definition}
\theoremstyle{remark}

\def\Alinea#1{\hfill\break%
  \hbox to \parindent{\hss{\textup{#1}}\enspace}\ignorespaces}
\def\alinea#1{\noindent%
  \hbox to \parindent{\hss{\textup{#1}}\enspace}\ignorespaces}

\def\Chi{\setbox0=\hbox{$\chi$} \mathord{\raise\dp0\hbox{$\chi$}}}

\def\H{\mathbf{H}}

\def\R{\mathbf{R}}
\def\S{\mathbf{S}}

\def\RP^#1{\mathbf{P}^{#1}(\mathbf{R})}

\def\qed{\hbox{\ Q.E.D.}
\vspace{.3cm}}

\def\diff{\mathrm{Diff}}
\def\sup{\mathrm{sup}}
\def\supp#1{\mathrm{supp(#1)}}
\def\inf{\mathrm{inf}}
\def\Le{\mathcal{L}}
\def\Ma{\mathcal{M}}
\def\Mo{\mathbf{M}}
\def\H{\mathcal{H}}
\def\Hq{\hat\H}
\def\Cm{\mathcal{C}^-}
\def\Cp{\mathcal{C}^+}
\def\mr#1{\;\mathrm{#1}\;}
\def\I{\mathcal{I}}

\def\adresse#1{\def\Y{\egroup\hbox\bgroup\sl} \par \noindent
  \hbox to \textwidth {\hfill \vbox {\small \hbox \bgroup #1 \egroup}}}
\title{Invariance of global solutions of the Hamilton-Jacobi equation}
\author{Ezequiel \textsc{Maderna}
}
\date{(in {\sl Bull. Soc. Math. France} (130) n.4, 2002)}

\begin{document}

\maketitle

\begin{abstract}
We show that every global viscosity solution of the Hamilton-Jacobi equation
associated with a convex and superlinear Hamiltonian on the cotangent bundle of
a closed manifold is necessarily invariant under the identity component of the
group of symmetries of the Hamiltonian (We prove that this group is a compact
Lie group). In particular, every Lagrangian section
invariant under the Hamiltonian flow is also invariant under this group.
\end{abstract}

\section{Introduction}

Let $M$ be a closed manifold, and let $H:T^*M\to\R$ be a C$^\infty$
Hamiltonian that is $C^2$-strictly convex and superlinear on the fibers of the cotangent
bundle $\pi^*:T^*M\to M$.

In \cite{F1}, generalizing work by Lions, Papanicolau \& Varadhan,
Fathi proved the existence of global viscosity solutions, also
called weak KAM solutions, of the Hamilton-Jacobi equation
$$H(x,d_xu)=c$$
and that these solutions only exist for the value $c=c(L)$, which equals
Ma\~{n}\'{e}'s critical value of the associated Lagrangian. The latter can
also be characterized in terms of Mather's minimizing measures
(see \cite{Ma}\cite{M1}). The solutions are given modulo constants
by the fixed points of the Lax-Oleinik semigroups
$T^-_t$ and $T^+_t$ (see below for the definition of $T^-_t$ and $T^+_t$).
Now let $\mathcal{S}_-$ and $\mathcal{S}_+$ be the
set of weak KAM solutions of $T^-_t$ and $T^+_t$ respectively. One has that
$\mathcal{S}_-\cap\mathcal{S}_+=\mathcal{S}$, the set of classical
solutions of the Hamilton-Jacobi equation, i.e. of class C$^1$.

The weak KAM solutions are very useful in the study of the dynamics of the
Hamiltonian vector field $X_H$ associated with $H$ (see also \cite{CDI}
\cite{F2}).

We will denote by $\Gamma_H$ the group of diffeomorphisms of $M$ of
class C$^1$ that preserve $H$, more precisely
$$\Gamma_H=\{g \in \diff^1(M)\,/\,
H(g(x),p)=H(x,p\circ d_xg) \; \forall x \in M, p \in T_{g(x)}^*M \}$$
endowed with the topology of uniform convergence. Let $\Gamma_H^0$
be the identity component of $\Gamma_H$. We shall prove in section
\ref{grupo} that $\Gamma_H$ is a compact Lie group.

In \cite{FM}, the proof of the existence of weak KAM solutions is
generalized to the case when $M$ is not necesarily compact, with
the additional hypothesis of uniform superlinearity (with respect
to a complete Riemannian metric) of the Hamiltonian and its
associated Lagrangian. In \cite{FM} we also show the existence of
$\Gamma_H$-invariant weak KAM solutions for values of the constant
greater or equal than a certain value
$c_{inv}\geq c(L)$. It follows that if $M$ is compact $c_{inv}=c(L)$.
We will prove later these facts in a slightly simplified way
using compactness.

On the other hand, if $M$ is not compact the inequality $c_{inv}\geq c(L)$
could be strict. This follows from the examples given by G.\& M. Paternain
\cite{PP} on the universal cover of closed surface of genus $2$.

In this paper we show:

\begin{theorem} Let $M$ be a closed manifold, and let
$H:T^*M\to\R$ be a C$^\infty$
Hamiltonian that is convex and superlinear on the fibers of the cotangent
bundle of $M$. If $\Gamma_H$ is the symmetry group of $H$,
then every weak KAM solution of $H$ is $\Gamma_H^0$-invariant,
where $\Gamma_H^0$ denote the identity component of $\Gamma_H$.
\end{theorem}

In general, Hamiltonians have trivial symmetry groups like general
Riemannian metrics which usually have trivial isometry groups. However we find
Hamiltonian systems with symmetries quite often in the applications,
and these symmetries are very useful for a detailed study of the
system. If the dimension of $\Gamma_H$ is sufficiently large we can
find all weak KAM solutions by integration, as we will see in the case
when $H$ is the Hamiltonian of the mechanical system determined by
the motion of a particle on the $n$-sphere $\S^n \subset \R^{n+1}$
under the effect of a potential
$U(x)=x_{n+1}$.
In this case we can reduce the problem to finding the weak KAM solutions
of the pendulum on the circle.

In particular, our result applies to the solutions of the
Hamilton-Jacobi equation of class C$^2$ that correspond to
exact Lagrangian sections of $T^*M$ invariant under the Hamiltonian flow
associated with $H$. To reduce the study of Lagrangian sections
to exact ones, we shall recall in section \ref{promedios}
that given any cohomology class in $H^1(M,\R)$ there exist a
closed $\Gamma_H$-invariant $1$-form that represent the class.
Combining this result with theorem $1$ we obtain the following
corollary whose proof will also be given in section \ref{promedios}:

\begin{corollary}\label{corollary}
Every Lagrangian section of $T^*M$ invariant
under the Hamiltonian flow of $H$ is also invariant under $\Gamma_H^0$.
\end{corollary}


\section{Weak KAM solutions and Mather's set.}

Before giving the proof of Theorem 1, we briefly recall the
properties of the weak KAM solutions which we will use.
The details of the proofs can be found in \cite{F1} and \cite{F3}.

Let us introduce initially the Lagrangian corresponding to $H$ like
its convex dual on the tangent bundle of $M$:
$$L:TM \to \R, \;\; L(x,v) =\sup\{ p(v) - H(x, p): p \in T_x^*M \} $$
It is well-known that $L$ is also of C$^\infty$ class, strictly convex
and superlinear on the fibers, i.e. its second derivative
$\partial^2L/\partial v^2$ is definite positive everywhere and for all
$K \in \R$ there exists a constant
$C_K$ such that
$$\forall (x,v) \in TM, \;\;L(x, v)\geq K \Vert v \Vert + C_K. $$

The Legendre transform $\Le:TM \to T^*M$,
$$\Le(x, v)=(x, \frac{\partial L}{\partial v}(x, v))$$
is a diffeomorphism which conjugates the Euler-Lagrange flow defined by
$L$ on $M$, which is denoted $\phi^L_t$, with the
Hamiltonian flow of $H$.

The action of $L$ on a piecewise C$^1$ curve $\gamma:[a,b] \to M$ is as usual
$$A_L(\gamma)=\int_a^b L(\gamma(s),\dot\gamma(s))\,ds\;.$$
We will say that a function $u:M \to \R$
is dominated by $L+c$ (for certain value of $c \in \R$ and we will write
$u \prec L+c$) if for each piecewise C$^1$ curve
$\gamma:[a,b] \to M$ we have:
$$u(\gamma(b))-u(\gamma(a)) \leq A_L(\gamma) + c(b-a).$$

From the superlinearity of $L$ it is easy to deduce that dominated
functions are Lipschitz, with a Lipschitz constant which only depends,
once fixed the metric on $M$, on the constant $c$ and the Lagrangian;
in particular, in accordance with the Rademacher's theorem (see \cite{Zi}),
they are differentiable almost everywhere.

The main reason we are interested in dominated functions
is that they constitute a suitable space where Lax-Oleinik's
semigroups of operators can be studied. In this way, we will
obtain weak KAM solutions. Since, we already know that the solutions
are dominated by $L+c(L)$, where $c(L)$ is the critical value of $L$,
we can directly introduce the space
$$\H  = \{ u \in C^0(M,\R) : u \prec L+c(L)\}\;;$$
on this space we can define, for each $t \geq 0$, the non
linear operators
$$u \to T^-_tu \;,\;\; u \to T^+_tu \;,$$
$$T^-_tu(x) = \inf_{\gamma\in \Cm} \{u(\gamma(0))+ A_L(\gamma)\}$$
$$T^+_tu(x) = \sup_{\gamma\in \Cp} \{u(\gamma(t))- A_L(\gamma)\}$$
where
$$\Cm = \{\gamma: [0,t]\to M \mr{piecewise}C^1, \mr{with} \gamma(t)=x\}$$
$$\Cp = \{\gamma: [0,t]\to M \mr{piecewise}C^1, \mr{with} \gamma(0)=x\}\;.$$

From the definition, it follows the semigroup property
$T^-_t \circ T^-_s = T^-_{t+s}$ for all $t,s \geq 0$,
and that $T^-_t(u+c)=T^-_t(u)+c$ for all $c\in\R$.
On the other hand, it is clear that $u+c\in\H$ whenever $u\in\H$;
this allows us to define the quotient semigroup $\hat T^-_t$
acting on $\Hq$, the quotient set of $\H$ by the space of constant
functions, by
$$\hat T^-_t [u]=[T^-_tu]\,.$$
Analogously, we can define the quotient semigroup $\hat T^+_t$.

\begin{definition} {\bf (Weak KAM solution.)}
We say that $u\in\H$ is a global viscosity solution of the
Hamilton-Jacobi equation, or a weak KAM solution, if
$\hat T^-_t [u]=[u]$ or $\hat T^+_t [u]=[u]$ for all $t\in\R$.
We call $\mathcal{S}_-$ and $\mathcal{S}_+$ respectively the sets
defined by the above relations.
\end{definition}

The existence of these solutions is obtained in \cite{F1} through
the application of the fixed point theorem of Schauder and Tykhonov;
this requires to show the continuity of the semigroups and the
compactness of the convex $\Hq$.
In the same article, it is shown that the relations which define
$\mathcal{S}_-$ and $\mathcal{S}_+$ sets, are equivalent to
$$T^-_t u = u - c(L)t \;\;\mr{and}$$
$$T^+_t u = u + c(L)t \;\;\; \forall t \geq 0 \;,$$
and that weak KAM solutions verify the Hamilton-Jacobi equation
at every point where they are differentiable.

Weak KAM solutions are also characterized by the fact of being dominated by
$L+c(L)$ and by the existence of certain curves on which their variation
is maximal; more precisely,

\begin{proposition}\label{calibrantes}
(Fathi \cite{F1}) A function $u:M \to \R$ is in $\mathcal{S}_-$
if and only if:

a) $u \prec L+c(L)$, where $c(L)$ is the critical value of $L$,

b) For all $x \in M$ there exists an extremal of $L$,
$\gamma_x : (-\infty,0] \to M$\\
with $\gamma_x(0)=x$, and such that
$\forall t \geq 0$ we have,
$$u(x) - u(\gamma_x(-t))=
\int_{-t}^0 L(\gamma_x(s),\dot\gamma_x(s))\,ds + c(L)t \,.$$

Moreover, the set of differentiability points of a function
satisfying a) contains the points $x \in M$ for which there
exists $\epsilon > 0$ and an extremal
$\gamma:[-\epsilon, \epsilon ] \to M$,
such that $\gamma(0)=x$ and
$$u(\gamma(\epsilon)) - u(\gamma(-\epsilon)) = \int_{-\epsilon}^\epsilon
L(\gamma(s), \dot\gamma(s)) \, ds + 2\epsilon c(L) \,.$$
\end{proposition}

The characterization of the functions in $\mathcal{S}_+$ is analogous,
it is enough to replace the curves of the condition b) by curves
of the form $\gamma_x: [ 0, +\infty) \to M$ with $\gamma_x(0)=x$ along
which the equality is satisfied.
\vspace{.3cm}

Let now $\mu$ be a Borel measure on $TM$, invariant under the Euler-Lagrange
flow, and let $u:M \to \R$ be a $(L+c)$-dominated function.
Because of invariance of $\mu$ we have:
$$\int_{TM}(u\circ\pi\circ\phi_1^L - u\circ\pi )\,d\mu = 0$$
where $\pi:TM \to M$
is the canonical projection of the tangent bundle.
If one applies for each $v \in TM$, the domination of $u$ by $L+c$
to the curve $t \to \pi\circ\phi_s^L(v)$ with $t$ varying
in $[0,1]$, it results from it that
$$c + \int_{TM}\int_0^1 L\circ\phi_s^L \,ds\,d\mu \, \geq 0.$$
By reversing the order of integration and once again owing to the fact
that $\mu$ is invariant, it results that:
$$\int_{TM} L \, d\mu \, \geq -c \,.$$

In addition, if $u \prec L+c(L)$ is a solution in $\mathcal{S}_-$,
and $\gamma_x$ one of the extremals which is associated to him by
the proposition \ref{calibrantes}, one can build a $\phi_t^L$-invariant
probability measure on $TM$, supported on the $\alpha$-limit set of
$\gamma_x$; it is enough to take a weak limit when $t\to - \infty$
of the probability measures $\mu_t$ defined by
$$\int f \,d\mu_t = -\frac{1}{t}\int_t^0 f(\gamma_x(s))\,ds$$
for  $f:TM\to\R$ a continuous function.
It is easy to note that the measures $\mu$ thus built satisfy
$$\int_{TM} L \, d\mu = -c(L) \,.$$
This indeed shows well the next characterization of the critical value
of the Lagrangian:

\begin{enumerate}

\item{$c(L)$ is the least value of $c\in\R$ such that the set\\
$\{f:M\to\R / f \prec L+c \}$ is not empty.}

\item{$c(L) = -\inf \,\int_{TM} L \,d\mu$, where the infimum is
taken over all probability measures invariants by the Euler-Lagrange
flow.}

\end{enumerate}

According to Mather (reference \cite{M1}), a measure $\mu$ is said to be {\sl minimizing}
if it is Borel probability measure $\mu$, invariant by the Euler-Lagrange
flow, and $\int_{TM} L \,d\mu = -c(L)$. The Mather set $\tilde\Ma$
is the closure of the union of the supports of all minimizing measures;
it is thus compact, invariant by $\phi_t^L$, and it contains the
$\alpha$-limit sets of curves $\gamma_x$ referred to above.
If $(x,v)\in \tilde\Ma$, and $u$ be a solution in $\mathcal{S}_-$ or
$\mathcal{S}_+$, then for all $t,t' \in \R$ we have
$$u(\pi\phi_t^L(x,v)) - u(\pi\phi_{t'}^L(x,v)) =
\int_{t'}^t L(\phi_s^L(x,v))\,ds + c(L)(t-t').$$
This proves that $u$ is differentiable on $\Ma=\pi(\tilde\Ma)$
and that if $(x,v)\in \tilde\Ma$ then $d_xu=\Le(x,v)$ (\cite{F1},
proposition 3.); in particular $\pi:\tilde\Ma\to\Ma$ is one to one.


\section{Proof of Theorem 1.}

At first we show that $\Ma=\pi(\tilde\Ma)$ is invariant under $\Gamma_H$,
and that the restrictions of solutions in $\mathcal{S}_-$ to $\Ma$ are
$\Gamma_H^0$-invariant (the proof for the case of solutions in
$\mathcal{S}_+$ is analogous). Actually, it is evident that the
transformations which preserve $H$ also preserve $L$; to be more precise,
for all $g \in \Gamma_H$, and for all $(x,v) \in TM$ we have
$$L(g(x),d_xg(v))=L(x,v).$$
Then, $A_L(\gamma)=A_L(g\circ\gamma)$ for all piecewise C$^1$ curve
$\gamma:[a,b] \to M$. In particular extremal curves of $L$ are sent
under $\Gamma_H$ to extremal curves, showing that $dg:TM \to TM$ also preserves
the Euler-Lagrange flow.

If $\mu$ is an invariant measure, then $dg^*\mu$ is also invariant. Moreover, if
$\mu$ is minimizing $dg^*\mu$ is also minimizing since
$$\int_{TM}L\,d(dg^*\mu) = \int_{TM}L\circ dg \,d\mu = \int_{TM}L \,d\mu\,.$$
If $(x,v) \in \supp\mu$ for a certain minimizing measure $\mu$,
then it is clear that $(g(x),d_xg(v)) \in \supp{d(g^{-1})^*\mu}$.
This proves that $dg(\tilde\Ma)=\tilde\Ma$, because $dg$ is an homeomorphism.
Finally, one has that for all $g \in \Gamma_H$,
$$g(\Ma)=g(\pi(\tilde\Ma))=\pi(dg(\tilde\Ma))=\pi(\tilde\Ma)=\Ma \,.$$

Let us now consider two solutions $u$ and $u_0$ in $\mathcal{S}_-$ and note
$v=u-u_0$ their difference; by the last lines of the previous
section, both solutions are differentiable at every point $x$ in $\Ma$ and their
derivative must be equal to the Legendre transform of the unique vector
in $\tilde\Ma$ that projects on $x$. Then $v$ is differentiable in all
points of $\Ma$ and at these points $d_xv = d_xu - d_xu_0 = 0$.
By proposition \ref{liecompact}, the orbits by
$\Gamma_H^0$ of each point in $\Ma$ are closed submanifolds of $M$;
as $\Ma$ contains the orbits of its points, we have, if $x\in \Ma$
and $V_x=\Gamma_H^0(x)$ is the orbit of $x$, that $V_x \subset \Ma$,
and then $d(v\vert_{V_x})=0$.
As $V_x$ is a connected submanifold, $v$ must be constant on $V_x$.
That is to say, two arbitrary solutions differ by a constant on
the orbits contained in $\Ma$. But we know (lemma \ref{invariantsolutions},
see also \cite{FM} for a more general statement) that there exists
solutions $\mathcal{S}_-$ (or $\mathcal{S}_+$) invariant by $\Gamma_H$.
We therefore conclude that every weak KAM solution is invariant by
$\Gamma_H^0$ on $\Ma$.

We still have to show the invariance outside $\Ma$. Let us fix a solution
$u$ in $\mathcal{S}_-$, an element $g$ of $\Gamma_H^0$ and any point $x$
of $M$. We will show that $u(g(x))=u(x)$. Given $\epsilon >0$, we choose $\delta >0$
so that for all $x,x'\in M$ with $d(x,x')<\delta$ we have at the same time
$\vert u(x)-u(x')\vert < \epsilon$ and $\vert u(g(x))-u(g(x'))\vert < \epsilon$.
Let us now take the curve $\gamma_x$ given by proposition \ref{calibrantes},
and $y$ in its $\alpha$-limit set, which is contained in $\Ma$, whose existence we have
already seen; be as well $t<0$ so that $d(\gamma_x(t),y)<\delta$.
To simplify the notation we will write $\gamma=\gamma_x\vert_{[t,0]}$ and
$z=\gamma(t)$. We then have, through definition from $\gamma$ and
because of the domination of $u$ by $L+c(L)$, the following relations:
$$u(x)-u(z)=A_L(\gamma)-c(L)t \, , \;\; \textrm{et}$$
$$u(g(x))-u(g(z)) \leq A_L(g \circ \gamma) - c(L)t =
A_L(\gamma) - c(L)t\, ,$$
then $u(g(x))-u(x) \leq u(g(z))-u(z)$.
Moreover, we have $u(g(y))=u(y)$, when $y \in \Ma$. Since
$d(y,z) < \delta$, we have $\vert u(y)-u(z) \vert < \epsilon$ and
$\vert u(g(y))-u(g(z)) \vert <\epsilon$; this implies
$\vert u(g(z))-u(z) \vert < 2\epsilon$, hence $u(g(x))-u(x) \leq 2\epsilon$.
Considering that $\epsilon >0$ is arbitrary, we conclude that
$$u(g(x)) \leq u(x)$$
for all $x \in M$ and $g \in \Gamma_H^0$. If we applies this result to
$g^{-1}$, it is clear that we obtain the opposite inequality, which
proves the theorem.
\qed


\section{The group of transformations $\Gamma_H$}
\label{grupo}

Let us first observe that when $H$ is the Hamiltonian associated to
the metric on $M$, i.e. $H(x,p)=\frac{1}{2}\Vert p \Vert^2$, the
group of symmetries of $H$ is nothing but the group of the
isometries of $M$. A classic theorem by Myers and Steenrod
guarantees then that the group is in fact a Lie group (\cite{Ko}).
The aim of this section is to generalize this theorem to the group of
symmetries $\Gamma_H$ that we have defined. We recall that
$\Gamma_H$ consists of diffeomorphisms of
class C$^1$, and that it is endowed with the topology of
uniform convergence.

\begin{proposition} \label{liecompact}
The group of transformations $\Gamma_H$ is a compact Lie group.
\end{proposition}

This is an immediate consequence of the following lemma and of the
theorem due to Montgomery (\cite{Mo}, Th.2 p.208) which establishes that
all compact subgroup of the group of C$^1$ diffeomorphisms of
a manifold is a Lie group. The following lemma assure the
equicontinuity of the symmetries:

\begin{lemma}
There exists $K>0$ such that for all $g\in \Gamma_H$ and for all $x\in M$,
we have $\Vert d_xg \Vert \leq K$
\end{lemma}

{\sl Proof.} From the superlinearity of $H$ we know that there exists
a constant $C^*_1 \in \R$ so that for all $(x,p) \in T^*M$,
$$H(x,p) \geq \Vert p \Vert + C^*_1 \,.$$
If $g$ is such that $H(g(x),p)=H(x,p\circ d_xg)$, then
$$ H(g(x),p) \geq \Vert p\circ d_xg \Vert + C^*_1 \,.$$
On the other hand, the unitary cotangent bundle of $M$ is compact
as $M$ is compact. We deduct that for all $p \in T^*_{g(x)}M$ with
$\Vert p \Vert = 1$,
$$\Vert p \circ d_xg \Vert \leq A^*_1 - C^*_1 \;,$$
in which $A^*_1=\sup\{H(x,p)\,:\,(x,p)\in T^*M,\; \Vert p \Vert =1\,\}$.
This proves that $\Vert d_xg \Vert \leq K=A^*_1 - C^*_1$, a constant which does
not depend on $g$ nor on $x$, as we wanted to prove.
\qed

To show that the group is compact, we therefore only need to show
that it is closed in $C^0(M,M)$ and use the Ascoli's theorem.

\begin{lemma}
The group $\Gamma_H$ is closed within $C^0(M,M)$.
\end{lemma}

{\sl Proof.} Let $\delta>0$ be so that for all point $x\in M$, the
application $\pi\circ\phi^L_{\delta}$ defines a diffeomorphism of
the ball $B_x(0,1)=\{v\in T_xM\,:\Vert v \Vert\ < 1\,\}$ on its image;
we will denote by $\varphi_x$ this diffeomorphism and $U_x$ its image.
This choice of $\delta$ is possible since the Euler-Lagrange flow
comes from a second order equation on $M$ (this means
more precisely that $d\pi\circ X_L=id_{TM}$), and $M$ is compact.

Let $g$ be an element of $\Gamma_H$ and $x$ a point in $M$. We know that
$dg$ commutes with the flow $\phi^L_t$ for all $t\in \R$, as $g$ preserves
the extremals of $L$. That is to say, for all $(x,v)\in TM$, and for all
$t\in \R$ we have
$$\phi^L_t(d_xg(v))=d_{\pi\phi^L_t(v)}g\,(\phi^L_t(v))$$
so
$$\pi\circ\phi^L_t(d_xg(v))=g(\pi\circ\phi^L_t(v)) \,.$$
It follows that on a small enough neighbourhood $U\subset U_x$
of the point $x$ (which by the previous lemma we can clearly choose
independent of $g$),
$$g=\varphi_{g(x)}\circ d_xg \circ \varphi_x^{-1}\;.$$

Let $\{g_n\}\subset \Gamma_H$ be a sequence converging uniformely to $g:M\to M$.
If we fix $x\in M$, extracting a subsequence if necessary,,
we can suppose that $d_xg_n$ tends to the linear application
$\alpha \in L(T_xM,T_{g(x)}M)$. On the neighbourhood $U$
of $x$ the sequence
$$g_n=\varphi_{g_n(x)}\circ d_xg_n \circ \varphi_x^{-1}$$
converges uniformely towards both $g$ and
$$\varphi_{g(x)}\circ \alpha \circ \varphi_x^{-1}\,.$$
This shows that $g$ is C$^1$, and that $g_n \to g$ in the C$^1$ topology.
Since $H$ is continuous it follows that $g\in\Gamma_H$,
which proves the lemma.
\qed

Finally, let us remark that we have also proved that the elements of
$\Gamma_H$ are in fact of class C$^\infty$. It is not hard to see
that the proofs which have been given for theorem 1 and his
corollary are also valid for Hamiltonians of class C$^k$ with
$k\geq 3$. In this case, both the Euler-Lagrange flow and the elements
of the group $\Gamma_H$ are of class C$^{k-1}$. In the proof of
the last lemma, the fact that the Hamiltonian is at least of class C$^3$
is used to guarantee that the maps $\pi\circ \phi^L_\delta$
are indeed diffeomorphisms.


\section{Invariant means}
\label{promedios}

We will now prove the existence of invariant solutions, a fact
that we used in the proof of theorem 1. As we have already said,
even if in \cite{FM} the existence for Hamiltonians on not necessarily
compact manifolds is established, we include here a proof in the
compact case.
Essentially, we shall prove that the vector space formed by
$\Gamma_H$-invariant functions is stable under the Lax-Oleinik
semigroups, and that they must have a fixed point in the intersection
of this space with the convex set of functions dominated
by $L+c(L)$ which is also stable under the Lax-Oleinik semigroups.
To see that this intersection is not empty, we have to
take invariant means of dominated functions.
This same averaging technique, proves as is well known that each cohomology
class contains invariant closed forms; once this will be done, we will be
able to give the proof of corollary \ref{corollary}.

We denote by $\Mo$ the invariant mean on the space of continuous functions
defined on $\Gamma_H$ induced by the normalized Haar integral:
$$\Mo : C^0(\Gamma_H,\R) \to \R \;,$$
$$\Mo(\psi) = \int_{\Gamma_H} \psi(g)\,d\lambda(g) \;,$$
where $\lambda$ is the left invariant Haar measure on $\Gamma_H$.

If we take a continuous function on $M$ and we average with $\Mo$ their
restrictions to each orbit of the $\Gamma_H$-action, it is clear that
we obtain an invariant continuous function on $M$; we shall also
denote by $\Mo$ the operator on $C^0(M,\R)$ thus defined.
That is, for a fixed function $u\in C^0(M,\R)$, and $x\in M$,
$$\Mo u (x)= \int_{\Gamma_H} u(g(x))\,d\lambda(g) \;.$$

\begin{lemma}
For each $c \in \R$, the convex set of functions dominated by $L+c$
is stable under the invariant mean operator $\Mo$.
\end{lemma}

{\sl Proof.}
Fix $c \in \R$, and take any function $u \prec L+c$.
Let $\gamma: [a,b] \to M$ be a piecewise C$^1$ curve.
If we apply the domination of $u$ to a translated of
$\gamma$ by a symmetry $g$ in $\Gamma_H$, we have
$$u(g(\gamma(b)))-u(g(\gamma(a))) \leq A_L(g \circ \gamma) + c(b-a).$$
But $A_L(g \circ \gamma)=A_L(\gamma)$ for all $g$ as the action
is preserved by symmetries. Meanning these inequalities we obtain
\begin{eqnarray*}
\Mo u (\gamma(b))-\Mo u (\gamma(a)) &=&
\int_{\Gamma_H} u(g(\gamma(b)))-u(g(\gamma(a)))\,d\lambda\\
&\leq & A_L(\gamma) + c(b-a),
\end{eqnarray*}
which proves well that $\Mo u \prec L+c$.
\qed

\begin{lemma}\label{invariantsolutions}
The sets of weak KAM solutions $\mathcal{S}_-$ and $\mathcal{S}_+$ both contain
$\Gamma_H$-invariant functions.
\end{lemma}

{\sl Proof.} Let us remember that we have denoted $\H$ the convex subset of
$C^0(M,\R)$ formed by the $L+c(L)$ dominated functions; let us name
$\I$ the linear space of $\Gamma_H$-invariant functions, and we define
$$\H_{inv} = \H \cap \I \;.$$
By the previous lemma, this intersection is not empty,
since $\H$ is itself not empty. If we take the quotient of $\H_{inv}$ by the space
of constant functions, we obtain a non empty convex set $\Hq_{inv}\subset\Hq$.
Since $\I$ is closed in $C^0(M,\R)$ and $\Hq$ is compact and convex, we obtain
that $\Hq_{inv}$ is also a compact and non-empty convex subset of $\hat C^0(M,\R)$, the
quotient of the space of continuous functions by the subspace of constant functions.

Moreover, the Lax-Oleinik semigroups preserve the space of invariant functions:
if $t\geq 0$ and $u\in \I$, then once $x\in M$ and $g \in \Gamma_H$ are fixed, we know
that for all piecewise C$^1$ curve $\gamma:[0,t] \to M$ with $\gamma(t)=x$,
$$u(g \circ \gamma (0)) + A_L(g \circ \gamma) = u(\gamma(0)) + A_L(\gamma) \;.$$
Taking the infimum over all these curves, it results that $T^-_tu(g(x)) \leq T^-_tu(x)$,
and consequently $T^-_tu(x)=T^-_tu(g^{-1}g(x)) \leq T^-_tu(g(x))$. We have then
proved that $T^-_tu \in \I$.

Therefore, the quotient semigroup $\hat T^-_t$ leaves $\Hq_{inv}$ invariant.
Its continuity is deduced from the one of $T^-_t$.
Applying the theorem of Schauder and Tykhonov we obtain a common fixed point in $\Hq_{inv}$
to the whole semi-group, that is to say, a class $\hat u\in \Hq_{inv}$ such that
$\hat T^-_t \hat u=\hat u$ for all $t\in\R$. It is clear that any function in this
class is an invariant weak KAM solution.
\qed

\begin{lemma}
Let $\Omega \in H^1(M,\R)$ be a cohomology class. There exists a $1$-form
$\omega_0$ in $\Omega$ which is invariant by the action of $\Gamma_H^0$.
\end{lemma}

{\sl Proof.}
The connectedness of $\Gamma_H^0$ guarantees that all its elements are isotopic
to the identity map on $M$; consequently, for all $g\in\Gamma_H^0$ and all
closed form $\omega$, we have that the pull-back of $\omega$ by $g$, which
was noted $g^*\omega$, is homotopic to $\omega$. In particular, they must
be cohomologous. To say this in a more suitable way, the affine
action of $\Gamma_H^0$ on the space of closed forms preserves the cohomology
classes. As $\Gamma_H^0$ is compact, and cohomology classes are also
affine spaces, each one of these cohomology classes must have at least one
point fixed under $\Gamma_H^0$, i.e. an invariant form.
For $1$-form $\omega \in \Omega$
we can give the fixed point explicitly, using the
Haar measure, this time normalized on $\Gamma_H^0$:
$$\omega_0=\frac{1}{\lambda(\Gamma_H^0)}
\int_{\Gamma_H^0}g^*\omega \,d\lambda(g) \;.$$
The form $\omega_0$ is in $\Omega$ and clearly it verifies
$g^*\omega_0=\omega_0$ for all $g \in \Gamma_H^0$.
\qed

{\bf Proof of corollary \ref{corollary}.}
\vspace{.3cm}

The graph of a C$^1$ section $\omega :M\to T^*M$ is a Lagrangian
submanifold for the standard symplectic form of $T^*M$ if and only if
$\omega$ is a closed $1$-form on $M$. Moreover, this
submanifold is invariant by the flow $\phi^H_t$, if and only if $H\circ \omega$
is constant.
If $\Omega = [\omega] \in H^1(M,\R)$ is the cohomology class of $\omega$, and
$\omega_0$ the $1$-form in the same class given by the previous lemma, one has that
$\omega - \omega_0 = du$ for a certain differentiable function $u:M\to \R$.
The invariance of $\omega$ by the Hamiltonian flow of $H$ can then be written
$$\forall x \in M, \; H_0(x,d_xu)=c$$
where we have defined $H_0(x,p)=H(x,p+ \omega_0(x))$. It is easy to see that
$H_0$ is also convex and superlinear; moreover, $H_0$ define the same Hamiltonian
flow as $H$, since $d\omega_0=0$. So we can apply theorem 1 to obtain that $du$
is $\Gamma_{H_0}^0$- invariant; but $\Gamma_{H_0}^0=\Gamma_H^0$ as $\omega_0$ is
$\Gamma_H^0$-invariant.
\qed


\section{The Spherical pendulum.}

We give here an example where the group of symetries of the system
is large. In that case, using the results obtained,
we find all the weak KAM solutions by integration.
Let us begin by considering
the $n$-dimensional sphere $\S^n$ naturally embedded in $\R^{n+1}$ and the
potential $U:\S^n \to \R$ which associates to each point its last
coordinate; that is to say,
$$U(x_1, \ldots , x_{n+1})= x_{n+1}\;.$$

We will study the weak KAM solutions of the mechanical Hamiltonian
which describes the motion of a punctual mass forced to move around
on the sphere under the action of the potential $U$. More precisely,
$$H(x,p)=\frac{1}{2}\Vert p \Vert ^2 + U(x)\;.$$
We observe that the symmetries of this Hamiltonian are the restrictions
to $\S^n$ of the orthogonal transformations of $\R^{n+1}$ which fix
the last coordinate. The group $\Gamma_H$ is then naturally identified
to the orthogonal group $O(n)$, and $\Gamma_H^0$ to the group $SO(n)$.
The associated Lagrangian is
$L(x,v)=\frac{1}{2}\Vert v \Vert ^2 - U(x)$. So, its critical value is
$c(L)=1$, and the unique minimizing measure is supported on the point
$(N,0)\in T\S^n$, where $N$ designates the point of the sphere where
$U$ achieves its maximum.
The orbits by the action of $\Gamma_H^0$, coincides with the potential
level sets, i.e. the $(n-1)$-dimensional spheres obtained as intersection
of $\S^n$ with the horizontal hyperplanes $x_{n+1}=k$.
As we know, weak KAM solution are determined by its values on the
Mather set, so in this case, as $\Ma=\{N\}$, we have a unique solution
in $\mathcal{S}_-$ modulo an additive constant. In order to determine
it, we shall call $u$ this solution and we shall assume that $u(N)=0$.

Let $\gamma_x:(-\infty,0]$ be an extremal curve of the Lagrangian, associated
to a certain point of $x \in \S^n$ by proposition \ref{calibrantes}.
The fact that $u$ is differentiable in $\gamma_x(t)$ with $t<0$
follows the same proposition. In particular, we must have
for all $t<0$, that $\dot\gamma_x(t)$ is orthogonal to the kernel of $d_xu$.
Now we can apply theorem 1, and deduce that the kernel of $d_xu$ is tangent to
the level sets of $U$. Therefore, the curve $\gamma_x$ must be contained
in the vertical plane generated by $x$ and $N$.
That means that $\gamma_x$ corresponds to a trajectory of the pendulum on the
circle determined by this plane and the sphere. We also have that
$\gamma_x(t) \to N$ when $t\to -\infty$, which proves the differentiability
of $u$ in the set $\S^n - \{-N\}$. On this set, the Hamilton-Jacobi
equation can be written
$$\Vert d_xu \Vert =\sqrt{2-2U(x)}\;,$$
which permits, reparametrizing $\gamma_x$ by its last coordinate, to
calculate $u$ explicitly by integration.

We conclude that the only weak KAM solutions are given by
$$u_{\pm}(x_1, \ldots , x_{n+1})=
u(N) \pm\int_{x_{n+1}}^1\,\sqrt{\frac{2-2s}{1-s^2}}\,ds \,.$$
Moreover, there are no differentiable solutions of the Hamilton-Jacobi
equation, and no Lagrangian section is preserved by the Hamiltonian flow
of this system.


\section{Other examples.}

Finally, it is convenient to observe that these results cannot be extended
to non compact manifolds, or to
Lagrangian submanifolds which are not necessarily graphs. To see this,
it is sufficient to consider the following examples; in both cases,
the considered Hamiltonian is the corresponding one to the Riemannian
metric, its flow the geodesic flow, and so, its critical value is $c(L)=0$.

In a non compact manifold, we can have global solutions of the
Hamilton-Jacobi equation for values of the constant greater than
the critical value.
In $\R^n$, the equation is written $\Vert d_xu \Vert =c$, and each affine
function $u:\R^n \to \R$ is a global solution. Except the constant
functions (i.e. the  solutions for the value of $c=0$), they are not
invariant under isometries. This is essentially the same example that
we could give on $\mathbf H^n$, the hyperbolic space, Busemann's functions
playing the role of affine functions.

On the other hand, corollary \ref{corollary} is false for general Lagrangian
submanifolds of $T^*M$. Let us consider two opposed points $x$ and $-x$ in
$\S^n \subset \R^{n+1}$, and let us define the functions $u_+=d(\cdot, x)$
and $u_-=d(\cdot, -x)$,
where $d$ is the Riemannian distance on $\S^n$. If we call $U\subset\S^n$
the complementary open set of $\{x,-x\}$, we have that these two solutions
are differentiable on $U$, therefore their derivatives define two
Lagrangian graphs in $T^*U$, that we shall denote $G_+$ and $G_-$, both
diffeomorphic to the product $\S^{n-1}\times (0,1)$. Let also
$S_+=\{p\in T^*_x\S^n, \; \Vert p \Vert =1\}$ and
$S_-=\{p\in T^*_{-x}\S^n, \; \Vert p \Vert =1\}$.

We define now $N\subset T^*\S^n$ as the union of these two graphs,
and the two spheres $S_+$ and $S_-$. If we observe that $\Le^{-1}(N) \subset T\S^n$
is the set of all unitary tangent vectors which define geodesics passing by
$x$ and $-x$, we see that $N$ is an embedded submanifold of $T^*\S^n$,
diffeomorphic to the product $\S^{n-1}\times \S^1$, and preserved by the
Hamiltonian flow.
We observe that for $t=\frac{\pi}{2}$,
the flow $\phi^H_t$ sends $G_+$ and $G_-$ into open neighborhoods in $N$
of the spheres $S_+$ and $S_-$.
As the Hamiltonian flow preserves the symplectic structure of the
cotangent bundle, we have proved that $N$ is Lagrangian.
Clearly, $N$ is not invariant under the identity component of the
isometry group of $\S^n$.


\end{document}